\documentclass[10pt]{amsart}
\usepackage{amssymb, enumitem}
\usepackage[all]{xy}
\usepackage{aliascnt}
\usepackage{mathtools}
\usepackage[colorlinks=true, linkcolor=blue]{hyperref}
\usepackage[english]{babel}
\usepackage{tikz, tikz-cd}
\usetikzlibrary{matrix,decorations.pathmorphing,arrows}
\tikzset{commutative diagrams/.cd,arrow style=tikz,diagrams={>=stealth'}}
\newcounter{TmpEnumi}

\newcommand{\Z}{{\mathbb{Z}}}
\newcommand{\R}{{\mathbb{R}}}

\newcommand{\N}{{\mathbb{N}}}

\newcommand{\Aut}{{\mathrm{Aut}}}

\theoremstyle{definition}
\newtheorem{lma}{Lemma}[section]

\numberwithin{equation}{section}

\newaliascnt{thmCt}{lma}
\newtheorem{thm}[thmCt]{Theorem}
\aliascntresetthe{thmCt}

\newaliascnt{corCt}{lma}
\newtheorem{cor}[corCt]{Corollary}
\aliascntresetthe{corCt}

\newaliascnt{propCt}{lma}
\newtheorem{prop}[propCt]{Proposition}
\aliascntresetthe{propCt}

\newtheorem*{thm*}{Theorem}
\newtheorem*{qst*}{Question}
\newtheorem*{cor*}{Corollary}
\newtheorem*{prop*}{Proposition}

\newtheorem{thmx}{Theorem}

\newaliascnt{pgrCt}{lma}

\aliascntresetthe{pgrCt}

\newaliascnt{dfCt}{lma}
\newtheorem{df}[dfCt]{Definition}
\aliascntresetthe{dfCt}

\newaliascnt{remCt}{lma}
\newtheorem{rem}[remCt]{Remark}
\aliascntresetthe{remCt}

\newaliascnt{remsCt}{lma}

\aliascntresetthe{remsCt}

\newaliascnt{egCt}{lma}

\aliascntresetthe{egCt}

\newaliascnt{egsCt}{lma}

\aliascntresetthe{egsCt}

\newaliascnt{qstCt}{lma}

\aliascntresetthe{qstCt}

\newaliascnt{pbmCt}{lma}

\aliascntresetthe{pbmCt}

\newaliascnt{notaCt}{lma}
\newtheorem{nota}[notaCt]{Notation}
\aliascntresetthe{notaCt}

\AtBeginDocument{   \def\MR#1{}
}

\title{Stable and real rank for crossed products by finite groups}
\author{Parisa Elyasi and Nasser Golestani}

\address{Department of Pure Mathematics, Faculty of Mathematical Sciences,
	Tarbiat Modares University,
	Tehran\\ Iran}
\email{p.elyasi@modares.ac.ir}

\address{Department of Pure Mathematics, Faculty of Mathematical Sciences,
	Tarbiat Modares University,
	Tehran\\ Iran}
\email{n.golestani@modares.ac.ir}


\begin{document}
	\maketitle
	\begin{abstract}
		A long-standing open question in the theory of group actions on C*-algebras is the stable rank of the crossed product. Specifically, N.~C.~Phillips
		asked that if a finite group $G$ acts on a simple unital C*-algebra $A$ with stable rank one, does the crossed product have stable rank one?
		A similar question can be asked about the real rank. Most of the existing partial answers contain
		a reasonable restriction (mainly, a Rokhlin-type property) on the action and assumptions on $A$. We remove all extra assumptions
		on $A$ (for instance, stable finiteness
		and that the order on projections over $A$ is determined by traces)
		and we prove that if the action has the tracial Rokhlin property
		and $A$ is simple and $\sigma$-unital with stable rank one or real rank zero,
		then so do the crossed product and the fixed point algebra.
		Moreover, we show that if the Kirchberg's central sequence algebra $\mathrm{F}(A)$ has
		real rank zero, then the weak tracial Rokhlin property is equivalent to the tracial Rokhlin
		property for actions on simple unital separable  C*-algebras $A$.
		
	\end{abstract}
	\section{Introduction}
	
	
	
	The study of group actions on C*-algebras and the structure of crossed products has always been a significant research area in operator algebras.
	A long-standing open question in this area is the stable rank of the crossed product.
	Blackadar asked in \cite[Question~8.2.3]{Bla90} that if
	$\alpha: G \to \Aut (A)$ is an action of a finite group $G$  on an AF~algebra $A$,
	does the crossed product $A\rtimes_{\alpha} G$ have stable rank one?
	The answer to this question is not known even if $G=\Z_2$. Pointwise outerness is not good enough for the kind of structural results for crossed products by finite groups. Example~8.2.1 of \cite{Bla90} gives an example of a pointwise outer action $\alpha$ of $\Z_2$ on a separable unital
	(nonsimple) C*-algebra $A$ such that $A$ has stable rank one but $A\rtimes_{\alpha} G$ has stable rank two.
	For actions on simple C*-algebras, Phillips
	asked in \cite[Problem~10.2.6]{GKPT} that if a finite group $G$  acts on a simple unital C*-algebra $A$ with stable rank one,
	does the crossed product $A\rtimes_{\alpha} G$ have stable rank one?
	A similar question can be asked about the real rank.
	Example~9 of \cite{Ell93} contains a pointwise outer action $\alpha$ of $\Z_2$ on a simple
	unital $\mathrm{AF}$ algebra $A$ such that $A\rtimes_{\alpha} G$ does not have real rank zero.
	
	\medskip
	

	Positive answers to these questions can be obtained if one puts reasonable restrictions on the action, for instance,
	the Rokhlin property \cite{OP12, sa15}. This property imposes K-theoretical obstructions
	on the underlying algebra.
	The tracial Rokhlin property  for finite group actions  on simple unital C*-algebras, suggested by   Phillips  \cite{Ph11},
	is generic
	in which the Rokhlin projections are not required to add up to unity,
	but rather just up to a small error in trace.
	Related to the classification program, the tracial Rokhlin property is useful to answer the permanence problem whether the crossed product share the same property
	as the original algebra  \cite{osphill06}.
	
	
	In \cite[Theorem~3.4]{Osaka} it is shown that if $\alpha$ is an action of a finite group $G$
	on a simple separable unital C*-algebra $A$ such that the order on projections over $A$ is determined by traces
	and $\alpha$ has the tracial Rokhlin property, then $\mathrm{sr}(A)=1$ implies
	$\mathrm{sr}(A\rtimes_{\alpha} G)=1$.
	A similar result \cite[Corllary~3.6]{Osaka} holds for real rank zero if moreover $A$ is ``stably finite."
	Note that, though the assumption of ``order on projections over $A$ is determined by traces"
	is not mentioned explicitly in the statement of these two results of \cite{Osaka},
	it is necessary for their proof since the authors use their Theorem~3.3. Also, in  \cite{FanFang}, these two results are claimed
	to be true without these assumptions, but we could not follow their proof.
	
	\medskip
	
	We remove all these extra assumptions and prove the following theorem which also covers the nonunital case
	(Theorems~\ref{thm_rr0} and \ref{sr1}):
	
	\begin{thmx}\label{12sr}
		\label{thma}	Let $A$ be a simple $\sigma$-unital C*-algebra, let $G$ be finite group, and let
		$\alpha:G\rightarrow\text{Aut}(A)$ have the tracial Rokhlin property.
		\begin{enumerate}
			\item If $A$ has stable rank one, then  $A\rtimes_{\alpha} G$ and  $A^{\alpha}$ have stable rank one.
			\item If $A$ has real rank zero, then $A\rtimes_{\alpha} G$ and $A^{\alpha}$ have real rank zero.
		\end{enumerate}
	\end{thmx}
	
	To achieve this, using the machinery of the Cuntz subequivalence
	(which, in particular, enables us to avoid using traces and to cover the
	nonunital case), we first generalize
	Lemma~3.1 of \cite{Archey} which is a tracial version of
	Theorem~3.2 of \cite{OP12} determines the local structure of the crossed product.
	In particular, we replace the assumption of real rank zero
	with Property~(SP), and we   remove assumptions of stable finiteness and that the order on projections over
	$A$ is determined by traces (see Proposition~\ref{prop_arc}).
	Then, by this and using the notion of  tracial  approximation  introduced by Elliott and Niu in
	\cite{EN08},  we obtain an analogue of Theorem~3.3 of \cite{Osaka} for simple unital C*-algebras
	but with no extra hypothesis (see Theorem~\ref{cor_ot}).  The final step
	is to pass from the unital to nonunital case which is provided in Proposition~\ref{propcross}.
	\medskip

	The tracial Rokhlin property for actions requires that underlying C*-algebras have rather enough projections.
	Instead, the weak tracial Rokhlin property uses positive contractions rather than projections.
	We do not know whether Theorem~\ref{12sr} holds for such actions. However,
	in some cases,   the tracial and weak tracial Rokhlin properties coincide. For instance, this is the case if the underlying algebra
	is simple with tracial rank zero \cite[Theorem~1.2]{forgol}.
	A natural question arises (see, e.g., \cite[Question~3.14]{Gar17}): if $A$ has sufficiently enough projections,
	does the weak tracial Rokhlin property imply the tracial Rokhlin property? In particular,
	one could ask if a simple C*-algebra $A$ has real rank zero and  $\alpha:G\rightarrow\text{Aut}(A)$
	has the weak tracial Rokhlin property, does it follow that $\alpha$ has the tracial Rokhlin property?
	Using the central sequence algebra, we give an answer to this question (Theorem~\ref{ultr}):

	\begin{thmx}\label{Ultrr0}
		Let $A$ be a simple separable C*-algebra and let $\alpha$ be an action
		of a finite group $G$ on $A$ with the weak tracial Rokhlin property. If $A'\cap A_{\omega}$ has
		real rank zero for some free ultrafilter  $\omega$ on $\N$, then $\alpha$ has the tracial Rokhlin property.
	\end{thmx}
	
	If $A$ is unital then $ A'\cap A_{\omega}=\mathrm{F}(A)$ is the Kirchberg's central sequence algebra \cite{Kirch}.
	As far as we know, the relation between the real rank of $A$ and that of $\mathrm{F}(A)$ is not known.
	If one could prove that $\mathrm{RR}(A)=0$ implies $\mathrm{RR}(\mathrm{F}(A))=0$, then
	combining this with Theorem~\ref{Ultrr0} gives a positive answer to the preceding question.
	A special case is that if $A$ is a simple unital separable purely infinite nuclear C*-algebra,
	then $\mathrm{F}(A)$ is simple and purely infinite, and hence $\mathrm{F}(A)$ has real rank zero,
	and so Theorem~\ref{Ultrr0} can be applied (Corollary~\ref{cor_Kirch}).
	
		
		

		\section{The Tracial Rokhlin Property}
		
		In this section   we   generalize Lemma~3.1 of \cite{Archey}
		which is a tracial version of
		Theorem~3.2 of \cite{OP12} (see Proposition~\ref{prop_arc} below).
		First, we fix some notation and then provide some preliminary lemmas.
		
		\begin{nota}
			Let  $A$ be a C*-algebra.
			\begin{enumerate}
				\item
				$A_{+}$ and $A_{sa}$ denote the positive cone
				and the set of self-adjoint elements of $A$, respectively. Also,
				$A^{+}$ denotes the  unitization of $A$
				(adding a  new identity  even if  $A$  is unital),
				while $A^{\sim} = A$ if $A$ is unital and
				$A^{\sim} = A^{+}$ if $A$ is nonunital.
				
				\item
				If $p $ and $q$ are projections in $A$, then we write
				$p\underset{{\mathrm{MvN}}}{\sim} q$ if $p$ is Murray-von Neumann equivalent to $q$.
				
				\item
				If $E$ and $F$ are  subsets of   $A$ and $\varepsilon>0$,
				then we write $E\subseteq_{\varepsilon} F$ if
				for every $a\in E$ there is $b\in F$ such that $\|a-b\|<\varepsilon$.
				Also, $a\in_{\varepsilon} F$ means $\{a\}\subseteq_{\varepsilon} F$.
				
				\item
				We write $\mathcal{K}=K(\ell^{2})$ and
				$M_{n}=M_{n}(\mathbb{C})$.
				
				\item
				For $a,b\in A_{+}$, we write $a\precsim_{A} b$
				if $a$ is Cuntz subequivalent to
				$b$, i.e., if there is a sequence $(v_n)_{n \in \mathbb{N}}$ in $A$ such that $\|a-v_{n}bv_{n}^{*}\|\to 0$.
				We write $a\sim_{A} b$ if both $a\precsim_{A} b$ and $b\precsim_{A} a$.
				Note that if $p$ is a projection then $p\precsim_{A} a$ if and only if there is a projection
				$q\in \overline{aAa}$ such that $p\underset{{\mathrm{MvN}}}{\sim} q$ \cite[Lemma~2.5]{forgol}.
				
			\end{enumerate}
		\end{nota}
		
		In  \cite{Ph11}, Phillips introduced the tracial Rokhlin property
		for finite group actions on simple separable unital C*-algebras. In our nonunital setting,
		we use the following generalization of the  (weak) tracial Rokhlin property
		to the nonunital case.
		
		\begin{df}[\cite{forgol}, Definition~3.1]\label{t.R.p}
			Let $\alpha \colon  G \rightarrow \mathrm{Aut}(A)$ be an action
			of a finite group $G$ on a simple C*-algebra $A$. Then
			$\alpha$ has the \emph{tracial Rokhlin property} if for  every
			finite subset $F \subseteq A$, every $\varepsilon > 0$, and
			every positive elements $x,y \in A$ with $\|x\|=1$,  there
			exists a family of mutually orthogonal projections $(p_{g})_{g \in G}$ in $A$ such that,
			with $p=\sum_{g \in G} p_{g}$, the following hold:
			\begin{enumerate}
				\item\label{t.R.p_it1}
				$\|p_{g}a-ap_{g}\| < \varepsilon$ for all $a\in F$ and all $g\in G$;
				
				\item\label{t.R.p_it2}
				$\|\alpha_{g}(p_{h})-p_{gh}\| < \varepsilon$ for all $g,h\in G$;
				
				\item\label{t.R.p_it3}
				$(y^{2}-ypy - \varepsilon)_{+} \precsim_{A} x$;
				
				\item\label{t.R.p_it4}
				$\|pxp\| > 1-\varepsilon$.
			\end{enumerate}
			The action $\alpha$ has the  \emph{weak tracial Rokhlin property} if we can arrange  $(p_{g})_{g \in G}$ above to be
			mutually orthogonal positive contractions.
		\end{df}
		
		\begin{rem}\label{alttrp}
			In Definition~\ref{t.R.p}, if instead of the assumption of mutual orthogonality of the projections $(p_g)_{g\in G}$,
			we assume that $(p_g)_{g\in G}$ are self-adjoint almost projections (that is,
			$\|p_g^2-p_g\|<\varepsilon$ for all $g\in G$) and they are almost orthogonal
			(that is, $\|p_gp_h\|<\varepsilon$ for all $g,h\in G$ with $g\ne h$), then we get an equivalent definition
			to the tracial Rokhlin property.
			This is because we can perturb these elements to be genuine orthogonal projections
			(see \cite[Lemma~2.5.5]{Lin.book} and \cite[Lemma~1.7]{Glimm60}).
		\end{rem}
		
		In the following lemma, we give a (seemingly) stronger equivalent definition of the
		tracial Rokhlin property. It says that we can arrange
		the projection $p$ in Definition~\ref{t.R.p} to be
		invariant, and we can take
		two different unknowns $x,z$ in
		Conditions~\eqref{t.R.p_it3} and
		\eqref{t.R.p_it4}  instead of  $x$. We will use this lemma in the proof of Proposition~\ref{prop_arc}.
		
		\begin{lma}\label{lemtrp1}
			Let $\alpha : G \rightarrow \mathrm{Aut}(A)$ be an action of a finite group $G$ on a simple
			C*-algebra $A$.
			Then  $\alpha$ has the tracial Rokhlin property if and only if
			the following holds. For
			every  $\varepsilon > 0$,  every finite subset
			$F \subseteq A$, and every positive elements $x,y,z \in A$ with $x\neq 0$ and $\|z\|=1$,
			there exists a family of orthogonal  projections $(p_{g})_{g \in G}$ in $A$
			such that, with $p=\sum_{g \in G} p_{g}$,
			we have $p\in A^{\alpha}$   and the following hold:
			\begin{enumerate}
				\item\label{lemtrp1_it1}
				$\|p_{g}a-ap_{g}\| < \varepsilon$ for all $a\in F$ and all $g\in G$;
				\item\label{lemtrp1_it2}
				$\|\alpha_{g}(p_{h})-p_{gh}\| < \varepsilon$ for all $g,h\in G$;
				\item\label{lemtrp1_it3}
				$(y^{2}-ypy - \varepsilon)_{+} \precsim_{A} x$;
				\item\label{lemtrp1_it4}
				$\|pzp\| > 1-\varepsilon$.
				\setcounter{TmpEnumi}{\value{enumi}}
			\end{enumerate}
		\end{lma}
		
		\begin{proof}
			The backward implication is obvious.
			For the forward implication, suppose that $\alpha$
			has the tracial Rokhlin property, and let $\varepsilon, F, x,y,z$ be as in the statement.
			We may assume that $F$ is contained in the closed unit ball of $A$. Let $n=|G|$ and put
			$\varepsilon_{0}=\min\{\frac{\varepsilon}{5}, \frac{1}{2}, \frac{\varepsilon}{1+2\|y\|^{2}}\}$.
			By \cite[Lemmas~2.5.1 and 2.5.4]{Lin.book}, there is $\delta_{0}>0$
			satisfying the following property: if  $b$ is a self-adjoint
			element in a C*-algebra $B$ and $e$ is a projection in $B$ such that $\|b-e\|<n\delta_{0}$,
			then there is a projection $p$ in $C^{*}(b)$ and a unitary $u\in B^{\sim}$ satisfying
			$ueu^{*}=p$ and $\|u-1\|<\varepsilon_{0}$, where $1$ denotes the unit of $B^{\sim}$.
			We may assume that $\delta_{0}<\varepsilon_{0}$.
			Choose $\delta$ with $0<\delta<1$ such that
			\[\left(\dfrac{\delta}{2-\delta}(1-\tfrac{\delta_{0}}{2})\right)^{2}>1-\delta_{0}.
			\]
			Put $z_{1}=(z^{ {1}/{2}}-\delta)_{+}$.
			Since $A$ is simple, \cite[Lemma~2.6]{Ph14} implies that there is a positive element
			$d\in \overline{z_{1}Az_{1}}$ such that $d\precsim x$ and $\|d\|=1$.
			Applying Definition~\ref{t.R.p} with $y$ and $F$ as given,
			with $ {\delta_{0}}/{2}$ in place of $\varepsilon$,
			and with $d$ in place of $x$, we obtain a family of orthogonal projections
			$(e_{g})_{g \in G}$ in $A$ such that, with $e=\sum_{g \in G} e_{g}$,
			the following hold:
			\begin{enumerate}
				\setcounter{enumi}{\value{TmpEnumi}}
				
				\item\label{lemtrp1_it5}
				$\|e_{g}a-ae_{g}\| < \frac{\delta_{0}}{2}$ for all $a\in F$ and all $g\in G$;
				
				\item\label{lemtrp1_it6}
				$\|\alpha_{g}(e_{h})-e_{gh}\| < \frac{\delta_{0}}{2}$ for all $g,h\in G$;
				
				\item\label{lemtrp1_it7}
				$(y^{2}-yey - \frac{\delta_{0}}{2})_{+} \precsim_{A} d$;
				
				\item\label{lemtrp1_it8}
				$\|ede\| > 1-\frac{\delta_{0}}{2}$.
				
				\setcounter{TmpEnumi}{\value{enumi}}
			\end{enumerate}
			We claim that
			\begin{enumerate}
				\setcounter{enumi}{\value{TmpEnumi}}
				
				\item\label{lemtrp1_it9}
				$\|eze\| > 1-\delta_{0}$.
				
				\setcounter{TmpEnumi}{\value{enumi}}
			\end{enumerate}
			To see this, first note that have
			$d\in \overline{z_{1}Az_{1}}\subseteq \overline{Az_{1}}=\overline{A(z^{{1}/{2}}-\delta)_{+}}$.
			Thus by \cite[Lemma~2.7]{forgol} there exists a sequence $(v_{n})_{n\in \mathbb{N}}$ in $A$
			such that $\|v_{n}z^{{1}/{2}}-d\|\to 0$ and
			$\|v_{n}\|\leq (\|d\|+\frac{1}{n})\delta^{-1}=(1+\frac{1}{n})\delta^{-1}$.
			Then $\|ev_{n}z^{\frac{1}{2}}e-ede\|\to 0$. Since $\|ede\|>1-\frac{\delta_0}{2}$ and
			$\delta<1$, there is $n\in \mathbb{N}$ such that
			$\|ev_{n}z^{\frac{1}{2}}e\|>1-\frac{\delta_0}{2}$ and $\frac{1}{n}<1-\delta$. Hence,
			\[
			1-\tfrac{\delta_0}{2}<\|ev_{n}z^{\frac{1}{2}}e\|\leq \|z^{\frac{1}{2}}e \| \|v_{n}\|
			\leq \|z^{\frac{1}{2}}e \| (1+\tfrac{1}{n})\delta^{-1}\leq
			\|z^{\frac{1}{2}}e \| (2-\delta)\delta^{-1}.
			\]
			Thus,
			\[
			\|eze\|=\|z^{\frac{1}{2}}e\|^{2}>
			\left(\dfrac{\delta}{2-\delta}(1-\tfrac{\delta_0}{2})\right)^{2}>1-\delta_0.
			\]
			This finishes the proof of the claim.
			
			Now we give an argument analogous to the proof of \cite[Lemma~1.17]{Ph11} to obtain
			the desired $p$ in $A^{\alpha}$. Put $b=\frac{1}{n}\sum_{g\in G}\alpha_{g}(e)$.
			Then by \eqref{lemtrp1_it6}, $\|b-e\|<n\delta_{0}$. Note that $b\in A^{\alpha}$.
			By the choice of $\delta_{0}$,
			there is a projection $p$ in $A^{\alpha}$ and a unitary $u\in A^{\sim}$ such that
			$ueu^{*}=p$ and $\|u-1\|<\varepsilon_{0}$ where $1$ denotes the unit of $A^{\sim}$.
			Put $p_{g}=ue_{g}u^{*}$, for all $g\in G$. Thus $(p_{g})_{g \in G}$ is a family of
			orthogonal  projections  in $A$ and $p=\sum_{g \in G} p_{g}$.
			We have
			\begin{enumerate}
				\setcounter{enumi}{\value{TmpEnumi}}
				
				\item\label{lemtrp1_it10}
				$\|p_{g}-e_{g}\| <2\varepsilon_{0}$, for all $g\in G$, and $\|p-e\| <2\varepsilon_{0}$.
				
				\setcounter{TmpEnumi}{\value{enumi}}
			\end{enumerate}
			Now we prove \eqref{lemtrp1_it1}--\eqref{lemtrp1_it4}.
			For \eqref{lemtrp1_it1}, by \eqref{lemtrp1_it5} and \eqref{lemtrp1_it10} we have
			\[
			\|p_{g}a-ap_{g}\|\leq \|p_{g}a-e_{g}a\|+\|e_{g}a-ae_{g}\|+\|ae_{g}-ap_{g}\|<5\varepsilon_{0}<\varepsilon.
			\]
			For \eqref{lemtrp1_it2}, using \eqref{lemtrp1_it6} and \eqref{lemtrp1_it10} we get
			\[
			\|\alpha_{g}(p_{h})-p_{gh}\|\leq \|\alpha_{g}(p_{h})-\alpha_{g}(e_{h})\|+
			\|\alpha_{g}(e_{h})-e_{gh}\|+\|e_{gh}-p_{gh}\|<5\varepsilon_{0}<\varepsilon.
			\]
			To see \eqref{lemtrp1_it3}, first by \eqref{lemtrp1_it10} we have
			\[
			\big \|(y^{2}-ypy)-(y^{2}-yey - \tfrac{\delta_{0}}{2})_{+}\big\| \leq \tfrac{\delta_{0}}{2}+
			\|ypy-yey\|<\varepsilon_{0}+2\varepsilon_{0}\|y\|^{2}\leq\varepsilon.
			\]
			Hence by \eqref{lemtrp1_it8} and \cite[Lemma~2.2]{KR02},
			$(y^{2}-ypy - \varepsilon)_{+}  \precsim_{A} (y^{2}-yey - \frac{\delta_{0}}{2})_{+} \precsim_{A} d \precsim_{A} x$.
			To prove \eqref{lemtrp1_it4}, by \eqref{lemtrp1_it9} and \eqref{lemtrp1_it10} we can compute
			\begin{align*}
				\|pzp\|&=\|eze+(p-e)zp+ez(p-e)\|\\
				&\geq \|eze\|-2\|p-e\|\\
				&>1-\delta_{0}-4\varepsilon_{0} \geq 1-\varepsilon.
			\end{align*}
			This finishes the proof.
		\end{proof}
		We need the following basic lemma and  Proposition~\ref{propemb} in the proof of Corollary~\ref{basic2}.
		
		\begin{lma} \label{basic}
			Let $A$ be a C*-algebra and $n \in \mathbb{N}$.
			Then for every self-adjoint  element $(a_{ij}) \in M_{n}(A)$ there exists
			$a \in A_{+}$ such that $(a_{ij}) \leq \mathrm{diag}(a,a,\ldots, a)$. More precisely, one can take
			$a=\sum_{i,j=1}^{n}|a_{ij}|$.
		\end{lma}
		
		\begin{proof}
			First, we observe that for any $b\in A$,
			\[
			\left|\begin{pmatrix}
				0 & b\\
				b^{*} & 0
			\end{pmatrix}\right|=
			\begin{pmatrix}
				|b^{*}| & 0\\
				0 & |b|
			\end{pmatrix}.
			\]
			Using this  we get
			\begin{equation}\label{equmat}
				\begin{pmatrix}
					0 & b\\
					b^{*} & 0
				\end{pmatrix}\leq
				\begin{pmatrix}
					|b^{*}| & 0\\
					0 & |b|
				\end{pmatrix}\leq
				\begin{pmatrix}
					|b|+|b^{*}| & 0\\
					0 & |b|+|b^{*}|
				\end{pmatrix}.
			\end{equation}
			Now let $B=(a_{ij})\in M_{n}(A)_{sa}$. For any $i,j$ with $1\leq i<j\leq n$  let
			$B_{ij}\in M_{n}(A)$ denote the matrix that has $a_{ij}$ and $a_{ji}$ in its $ij$-th and $ji$-th entries,
			respectively, and $0$ elsewhere. Also, let $C=\mathrm{diag}(a_{11},a_{22},\ldots, a_{nn})$.
			Put $a=\sum_{i,j=1}^{n}|a_{ij}|$ and
			$c=\sum_{i\neq j}|a_{ij}|$. Thus, by \eqref{equmat},
			$\sum_{i<j}B_{ij}\leq\mathrm{diag}(c,c,\ldots, c)$.
			We have
			\begin{align*}
				B=(a_{ij})=C+ \sum_{i< j}B_{ij} &\leq
				\mathrm{diag}(|a_{11}|,|a_{22}|,\ldots, |a_{nn}|)+\mathrm{diag}(c,c,\ldots, c)\\
				&\leq
				\mathrm{diag}(a,a,\ldots, a).\qedhere
			\end{align*}
		\end{proof}
		
		The following fact is known.
		
		\begin{prop}\label{propemb}(cf.~\cite[Corollary~4.1.6]{BO08})
			Let $\alpha$ be an action of a finite group $G$ on a C*-algebra $A$. Then there is an embedding
			of $A\rtimes_{\alpha} G$ into $M_{G}(A)$. More precisely, the map
			$\varphi: A\rtimes_{\alpha} G \to M_{G}(A)$ defined by
			\[
			\varphi\left(\sum_{s\in G} a_{s}s\right)=\left(\alpha_{s}^{-1}(a_{s	t^{-1}})\right)_{s,t\in G}
			\]
			is an injective $*$-homomorphism. If $A$ is unital then $\varphi$ is unital.
		\end{prop}

		\begin{cor}\label{basic2}
			Let $\alpha$ be an action of a finite group $G$ on a C*-algebra $A$.
			Then for every self-adjoint element $y \in A \rtimes _{\alpha}G$
			there exists $a \in A_{+}$ such that $y\leq a$ and  $y \in \overline{(A\rtimes _{\alpha}G)a}$.
		\end{cor}
		
		\begin{proof}
			Let
			$\varphi: A\rtimes_{\alpha} G \to M_{G}(A)$ be the
			$*$-homomorphism of
			Proposition~\ref{propemb}. Write $y=\sum_{s\in G} y_{s}s$.
			
			By Lemma ~\ref{basic},  $\varphi(y) \leq \mathrm{diag}(a,a,\ldots, a)$
			where $a= \sum_{s,t \in G}| \alpha_{s ^{-1}}(y_{st^{-1}})|$. For any $l\in G$ we have	
			\[
			\alpha_{l} (a)= \sum_{s,t \in G} |\alpha_{ls^{-1}}(y_{st^{-1}})|=
			\sum_{s,t \in G} |\alpha_{s^{-1}}(y_{{sl}t^{-1}})|=a,
			\]
			and  so $\varphi(a)=\mathrm{diag}(a,a,\ldots, a)$.
			Since $\varphi$ is a $*$-isomorphism,  we get $y \leq a$ in $A \rtimes _{\alpha}G$.
			For the last part, we have
			$y\in  \overline{a(A\rtimes _{\alpha}G)a}\subseteq \overline{(A\rtimes _{\alpha}G)a}$.
		\end{proof}
		
		We need the following two lemmas in the proof of Proposition \ref{prop_arc}. In the first lemma we extend a special case of \cite[Theorem~4.2]{JO98} to the nonunital case.
		\begin{lma} \label{proj.crossed}
			Let $A$  be a simple C*-algebra with Property~(SP) and let $\alpha$ be a pointwise outer action  of a
			discrete group $G$ on $A$.
			Then every nonzero hereditary C*-subalgebra of the reduced crossed product $A \rtimes _{\alpha, \mathrm{r}} G$
			has a nonzero projection which is Murray-von~Neumann
			equivalent (in $A \rtimes _{\alpha, \mathrm{r}} G$) to a projection of $A$.
			In particular, $A \rtimes _{\alpha, \mathrm{r}} G$ has Property~(SP).
		\end{lma}
		
		\begin{proof}
			Assume that $B$ is a hereditary C*-subalgebra of $A \rtimes _{\alpha, \mathrm{r}} G$.
			Since  $A$ is simple, by Proposition~1.4 of \cite{Pas}, the action $\alpha$ is pointwise spectrally nontrivial.
			Thus, Lemma~3.2 and Proposition~3.9 of \cite{Pas} imply that there exist a
			nonzero hereditary C*-subalgebra $E$
			of $A$ and an injective $*$-homomorphism $\phi \colon E \rightarrow B$ such that $\phi(x) \sim x$
			in $A \rtimes _{\alpha, \mathrm{r}} G$, for all $x \in E_{+}$.
			Since $A$ has Property~(SP), we can choose a nonzero projection $p$ in $E$.
			Hence, $\phi(p)$ is a projection in $B$ and $ p \precsim \phi(p)$ in $A \rtimes _{\alpha, \mathrm{r}} G$.
			Thus there is a projection $q\in  A \rtimes _{\alpha, \mathrm{r}} G$ such that
			$q \underset{{\mathrm{MvN}}}{\sim} p$ and $ q \leq \phi(p)$.
			Since $ \phi(p) \in B$ and $B$ is a
			hereditary C*-subalgebra of $A \rtimes _{\alpha, \mathrm{r}} G$, we have $q \in B$.
			This completes the proof.
			
			Alternatively, the statement follows from \cite[Theorem~4.2]{JO98} by taking $N=\{1\}$.
			Observe that the assumption that $A$ is unital is not used in the proof of \cite[Theorem~4.2]{JO98}.
		\end{proof}
		
		The following lemma should be known. The proof is similar to the proof of
		Lemma~2.5.3 of \cite{Lin.book}. We give the full proof for completeness.
		
		\begin{lma} \label{p.iso}
			For any $\varepsilon > 0$, there is a $\delta >0$ with the following property.  If $A$ is a C*-algebra,
			$p, q\in A$ are projections, and $x \in A$ satisfies
			\[
			\|x^{*}x-p\| < \delta, \ \
			\|xx^{*} -q\|< \delta,\ \ \|qxp-x\| < \delta,
			\]
			then there is a partial isometry $v  \in A$ such that
			$$v^{*}v=p, \ \ vv^{*}=q, \ \  \|x-v\| < \varepsilon.$$
		\end{lma}
		
		\begin{proof}
			There is $0 < \delta_{0} < 1$ such that if a real number $t $ satisfies $| t - 1 | < \delta_{0}$ then
			$|t^{-\frac{1}{2}}-1|<\frac{\varepsilon}{4}$. We set
			\[
			\delta=\min \left(\frac{\delta_{0}}{2\sqrt{2}+1},\frac{\varepsilon}{2}\right).
			\]
			
			Suppose that $p$, $q$, and $x$ are as in the statement. We put $y=qxp$. Then $qyp=y$ and
			$\|x-y\|=\|x-qxp\|<\delta$. Also, ${\|x\|^2\leq\|x^*x-p\|+\|p\|\leq\delta+1}$, and
			so $\|x\|\leq\sqrt{1+\delta}\leq\sqrt{2}$ and $\|y\|\leq\|x\|\leq\sqrt{2}$. Hence,
			${\|y^*y-x^*x\|<2\sqrt{2}\delta.}$
			Similarly, ${\|yy^*-xx^*\|<2\sqrt{2}\delta}$. Thus $\|y^*y-p\|<(2\sqrt{2}+1)\delta < 1$ and $\|yy^*-q\|<(2\sqrt{2}+1)\delta < 1$.
			Hence
			$y^*y$ is invertible in $pAp$ and $yy^*$ is invertible in $qAq$.
			Let $v=y|y|^{-1}$ where the inverse is taken in $pAp$. Then $v^*v=p$ and $vv^*=q$. We have
			\begin{align*}
				\|v-x\|&\leq\|v-y\|+\|y-x\|<\|y|y|^{-1}-yp\|+\delta
				<\sqrt{2}\, \||y|^{-1}-p\|+\delta.
			\end{align*}
			But $\|y^*y-p\|<(2\sqrt{2}+1)\delta < \delta_0$, and so by the choice of $\delta_0$
			we get $\||y|^{-1}-p\| = \|(y^*y)^{-\frac{1}{2}}-p\|<\frac{\varepsilon}{4}$.
			Therefore,
			$\|v-x\|<\frac{\sqrt{2}\varepsilon}{4}+\frac{\varepsilon}{2}<\varepsilon$,
			as desired.
		\end{proof}
		
		In the following result
		we generalize
		Lemma~3.1 of \cite{Archey} in some aspects.
		In particular, we replace the assumption of real rank zero
		with Property~(SP), and we   remove assumptions of
		stable finiteness and that the order on projections over
		$A$ is determined by traces.
		Also, we replace the unital assumption with $\sigma$-unital.
		If $A$ is unital, Property~(SP)
		is unnecessary (see Remark~\ref{prop_arc} below).
		Also, we do not need
		Condition~(6) of \cite[Lemma~3.1]{Archey}
		here. (In fact, a variant of this condition
		was part of the first definition of the tracial
		Rokhlin property for finite group actions,
		and then in the final definition, it was
		removed; see \cite[Remark~1.3]{Ph11}.)

		\begin{prop}\label{prop_arc}
			Let  $A$ be a simple C*-algebra with Property~(SP) and let $\alpha$ be an action of a finite group
			$G$ on $A$ with the tracial Rokhlin property. Then for every finite set
			$F \subseteq A\rtimes _{\alpha} G$, every $\varepsilon >0$, and every  $x,y,z \in (A\rtimes _{\alpha} G)_{+}$
			with $x\neq 0$
			and $\|z\|=1$,
			there exist a nonzero projection $ e \in A^{\alpha}$, a unital C*-algebra
			${D \subseteq e(A\rtimes_{\alpha}G)e}$, a nonzero projection $f \in A$, and an isomorphism
			${\phi: M_{n} \otimes fAf \rightarrow  D}$ where $n=|G|$, such that the following hold:
			\begin{enumerate}
				\item\label{prop_arc_it1}
				with $(e_{g,h})$ for $g, h \in G$ being a system of a matrix units for $M_{n}$,
				we have $\phi (e_{1,1} \otimes a)=a$ for all $a \in fAf$ and $\phi(e_{g,g}\otimes f) \in A$
				for all $g \in G$;
				\item\label{prop_arc_it2}
				with $(e_{g,h})_{g, h \in G}$ as above, we have
				$\|\phi(e_{g,g} \otimes a)- \alpha_{g}(a)\| \leq \varepsilon \|a\|$ for all $a \in fAf$
				and all $g \in G$;
				\item\label{prop_arc_it3}
				for any $a \in F$ we have $\| ea-ae\|<\varepsilon$ and
				there exist $b_{1}, b_{2} \in D$ such that $\|ea-b_{1}\| < \varepsilon$,
				$\|ae-b_{2}\| < \varepsilon$, and $\|b_{1}\|, \|b_{2}\| \leq \|a\|$;
				\item\label{prop_arc_it4}
				$\sum _{g \in G} \phi(e_{g,g} \otimes f)=e$;
				\item\label{prop_arc_it5}
				$ (y^{2}-yey-\varepsilon)_{+} \precsim x$
				in $A\rtimes _{\alpha}G$;
				\item\label{prop_arc_it6} $\|eze\|>1-\varepsilon$.
				\setcounter{TmpEnumi}{\value{enumi}}
			\end{enumerate}
		\end{prop}
		
		\begin{proof}
			Let $F$, $\varepsilon$, $x$, $y$, and $z$ be as in the statement. we may  assume that $\varepsilon<1$.
			We extend $\alpha:G\to \mathrm{Aut}(A)$
			to $\alpha:G\to \mathrm{Aut}(A^{\sim})$ by requiring that $\alpha_{g}(1)=1$ for all $g\in G$,
			where $1$ denotes the unit of $A^{\sim}$.
			Put $u_{g}=1 g\in A^{\sim}\rtimes_{\alpha}G$. Thus $u_{g}$ is a unitary,
			$\alpha_{g}(a)=u_{g}au_{g}^{*}$, and $a g=au_{g}$ for all $a\in A^{\sim}$ and all $g\in G$.
			For any $a\in F$ we can write
			$a = \sum_{g \in G} a_{g}u_{g}$. Let $S$ be a finite subset of $A$ containing all coefficients of elements
			of $F$.
			Let $\delta_{0}$ satisfy
			$0<\delta_{0}<1$,
			and choose $\delta_{1}$ according to
			Lemma~\ref{p.iso} for $\delta_{0}$ in place of $\varepsilon$.
			Choose $\delta>0$ such that
			\[
			\delta<\min\left\{\delta_{0},\delta_{1},\delta_{2}, \frac{\varepsilon}{10n^{3}}\right\}
			\]
			where $\delta_{2}$  is as follows. By Corollary~\ref{basic2}, there exists $y_1\in A_+$ such that ${y\in\overline{(A\rtimes _{\alpha}G)y_1}}$. Thus there is a nonzero element $w\in A\rtimes _{\alpha}G$
			such that
			${\|y-wy_{1}\|<\min\left\{1,\frac{\varepsilon}{4(2\|y\|+1)}\right\}}$. Put
			$\delta_{2}=\frac{\varepsilon}{2\|w\|^{2}}$.
			
			Choose $\eta$ with
			$\frac{1}{4}<\eta<1$ such that $\eta^{\frac{1}{4}}(1-\frac{\varepsilon}{2})>1-\varepsilon$.
			Consider the hereditary C*-subalgebra
			$\overline{(z-\eta)_{+}(A\rtimes _{\alpha}G)(z-\eta)_{+}}$ of $A\rtimes _{\alpha}G$.
			By Lemma~\ref{proj.crossed}, there is a nonzero projection $q\in A$ which is Murray-von~Neumann
			equivalent to a projection in $\overline{(z-\eta)_{+}(A\rtimes _{\alpha}G)(z-\eta)_{+}}$.
			Thus $q\precsim (z-\eta)_{+}$ in $A\rtimes _{\alpha}G$. Then by  \cite[Lemma~2.6]{forgol} there is an element
			$v\in A\rtimes _{\alpha}G$ such that
			$q=vzv^{*}$ and $\|v\|\leq \eta^{-\frac{1}{2}}$.
			Write $v=\sum_{g \in G} v_{g}u_{g}$. Put $V=\{v_{g}\mid g\in G\}$.
			
			Since $A$ is simple and has Property~(SP), by
			Lemma~\ref{proj.crossed} there exists
			a nonzero projection $q_{0}\in A$ which is Murray-von~Neumann equivalent
			(in $A\rtimes _{\alpha}G$) to a projection in
			$\overline{x (A\rtimes _{\alpha}G)x}$. Thus $q_{0}\precsim x$ in $A\rtimes _{\alpha}G$.
			Since $\alpha$ has the tracial Rokhlin property, applying Lemma~\ref{lemtrp1} with
			$\delta$ in place of $\varepsilon$, with $S \cup V$ in place of $F$, with $q_{0}$ in place of $x$, with
			$y_{1}$ in place of $y$, and with $q$ in place of $z$,
			we obtain a family of
			orthogonal  projections $(e_{g})_{g \in G}$ in $A$  such that with $e=\sum_{g \in G} e_{g}$
			we have $e\in A^{\alpha}$ and the following hold:
			
			\begin{enumerate}
				\setcounter{enumi}{\value{TmpEnumi}}
				\item\label{prop_arc_it7}
				$\|e_{g}a-ae_{g}\| < \delta$ for all $a\in S\cup V$ and all $g\in G$;
				\item\label{prop_arc_it8}
				$\|\alpha_{g}(e_{h})-e_{gh}\| < \delta$ for all $g,h\in G$;
				\item\label{prop_arc_it9}
				$(y_{1}^{2}-y_{1}ey_{1} - \delta)_{+} \precsim q_{0}$
				in $A\rtimes _{\alpha}G$;
				\item\label{prop_arc_it10}
				$\|eqe\| > 1-\delta$.
			\end{enumerate}
			Put $f=e_{1}$ and $w_{g,h}=u_{gh^{-1}}e_{h}$ for all $g,h\in G$.
			Note that by
			\eqref{prop_arc_it8} and \eqref{prop_arc_it10} we have $f\neq 0$.
			Also note that $w_{g,h}$ is in $A\rtimes _{\alpha}G$ since
			$A\rtimes _{\alpha}G$ is an ideal in $A^{\sim}\rtimes _{\alpha}G$.
			Then we have
			$\|w_{g,h}^{*}-w_{h,g}\|<\delta$ and $\| w_{g,h}w_{k,l}-\delta_{h,k}w_{g,l}\|<\delta$ for all
			$g,h,k,l\in G$. Also, $w_{1,1}=e_{1}\neq 0$. Similar to the proof of
			\cite[Lemma~3.1]{Archey} we have
			$\|w_{g,1}w_{g,1}^{*}-e_{g}\| < \delta$ and
			$\|w_{g,1}^{*}w_{g,1}-e_{1}\| < \delta$ for all $g\in G$.
			Since $\delta<\delta_{1}$, by the choice of $\delta_{1}$ there exist partial isometries
			$z_{g}\in A\rtimes _{\alpha}G$, $g\in G$,
			such that $\|z_{g}-w_{g,1}\|<\delta_{0}$,
			$z_{g}z_{g}^{*}=e_{g}$, and $z_{g}^{*}z_{g}=e_{1}$. We may take $z_{1}=e_{1}$.

			Let $(e_{g,h})_{g, h \in G}$ be a system of matrix units for $M_{n}$.
			Now define the map ${\phi: M_{n} \otimes fAf \rightarrow  A\rtimes _{\alpha}G}$ by
			$\phi(e_{g,h}\otimes a)=z_{g}az_{h}^{*}$, for all $g,h\in G$ and $a\in fAf$.
			Then $\phi$ is a $*$-homomorphism, $\phi(e_{g,g}\otimes f)=e_{g}\in A$, and
			$\phi(1\otimes f)=e$.
			Let $D$ be the range of $\phi$. Thus
			$\phi: M_{n} \otimes fAf \rightarrow  D$ is an isomorphism and $e=1_{D}$.
			Exactly similar to the proof of
			\cite[Lemma~3.1]{Archey} we get \eqref{prop_arc_it1}, \eqref{prop_arc_it2},
			\eqref{prop_arc_it4}, and the second part of
			\eqref{prop_arc_it3}.
			It remains to show \eqref{prop_arc_it5},
			\eqref{prop_arc_it6}, and the first part of
			\eqref{prop_arc_it3}. To see the first part of
			\eqref{prop_arc_it3} let
			$a\in F$ and write
			$a = \sum_{g \in G} a_{g}u_{g}$. So each $a_{g}$ is in $S$. By \eqref{prop_arc_it7} and that $e\in A^{\alpha}$ we get
			\begin{align*}
				\|ea-ae\|&= \Big\|\sum_{g\in G}ea_{g}u_{g}-\sum_{g\in G}a_{g}eu_{g}\Big\|\\
				&\leq
				\sum_{g\in G}\|ea_{g}-a_{g}e\|\\
				&\leq
				\sum_{g,h\in G}\|e_{h}a_{g}-a_{g}e_{h}\|\\
				&<n^{2}\delta<\varepsilon.
			\end{align*}
			The proof of \eqref{prop_arc_it5} is similar to the proof  \cite[Lemma~3.4]{forgol} and so is omitted.
			To prove  \eqref{prop_arc_it6}, first similar to
			the  inequalities above, we get
			$\|ev-ve\|<n^{2}\delta$. By this inequality,
			\eqref{prop_arc_it10}, and that
			$\|v\|\leq \eta^{-\frac{1}{2}}<2$ we obtain
			\begin{align*}
				1-\delta &<\|eqe\|=\|evzv^{*}e\|\\
				&\leq\|(ev-ve)zv^{*}e\|+
				\|vezv^{*}e\|\\
				&\leq \|ev-ve\|\cdot \|v\|+\|vezev^{*}\|+\|vez(ev^{*}-v^{*}e\|\\
				&\leq 2\|v\|\cdot \|ev-ve\|+\|v\|^{2}\|eze\|\\
				&<4n^{2}\delta+\eta^{-\frac{1}{4}}\|eze\|.
			\end{align*}
			Since $\delta\leq\delta_{0}<\frac{\varepsilon}{10n^3}<\frac{\varepsilon}{2(4n^{2}+1)}$, we get
			\[
			\|eze\|>\left(1-(4n^{2}+1)\delta\right)\eta^{\frac{1}{4}}>\left(1-\frac{\varepsilon}{2}\right)\eta^{\frac{1}{4}}
			>1-\varepsilon,
			\]
			which is \eqref{prop_arc_it6}. This finishes the proof.
		\end{proof}
		
		\begin{rem}\label{rmk_arc}
			In Proposition~\ref{prop_arc}, if moreover
			$A$ is unital then the assumption of Property~(SP)
			is unnecessary. Indeed, if $A$ does not have
			Property~(SP) then $\alpha$ has the Rokhlin property,
			by \cite[Lemma~1.13]{Ph11}. Now, the proof of
			\cite[Theorem~2.2]{Ph11} provides the desired
			$*$-isomorphism $\phi$ (called $\varphi$ in that
			proof). Also, note that
			if $A$ is unital then
			Condition~\eqref{prop_arc_it5} can be replaced with
			$ 1-e \precsim x$
			in $A\rtimes _{\alpha}G$, since we can take
			$y=1$ and use the fact that
			$ (1-e -\varepsilon)_{+}\sim 1-e$
			whenever $\varepsilon<1$.
		\end{rem}
		
		\section{Preservation of real rank zero and stable rank one}
		
		In this section, we provide the proof of Theorem~\ref{12sr}
		(see Theorems~\ref{thm_rr0} and \ref{sr1}).
		The following result is useful in the verification of passage of properties from a
		C*-algebra $A$ to crossed products $A\rtimes_{\alpha}G$.
		This is a generalization of Theorem~3.3 of \cite{Osaka}
		from which the assumptions of stable finiteness and that
		the order
		on projections over $A$ is determined by traces being removed.
		The proof is mainly based on Proposition~\ref{prop_arc}.
		
		\begin{thm}\label{cor_ot}
			Let $\mathcal{C}$ be a class of unital C*-algebras
			such that:
			\begin{enumerate}
				\item\label{cor_ot_it1}
				if $A \in \mathcal{C}$ and $B$ is a C*-algebra
				with $B\cong A$, then $B \in \mathcal{C}$;
				
				\item\label{cor_ot_it2}
				if $A \in \mathcal{C}$ and $n\in \mathbb{N}$
				then if $M_{n}(A) \in \mathcal{C}$;
				
				\item\label{cor_ot_it3}
				if $A \in \mathcal{C}$ and $p\in A$
				is a nonzero projection, then
				$pAp \in \mathcal{C}$.
				\setcounter{TmpEnumi}{\value{enumi}}
			\end{enumerate}
			Let $A$ be a simple unital C*-algebra with
			$A \in \mathrm{TA}\mathcal{C}$, that is, $A$ can be
			tracially approximated by the C*-algebras in
			$\mathcal{C}$ in the sense of \cite{EN08}. Let $\alpha$ be an action of
			a finite group $G$ on $A$ with the tracial
			Rokhlin property. Then
			$A\rtimes_{\alpha} G\in \mathrm{TA}\mathcal{C}$.
		\end{thm}
		
		\begin{proof}
			We need to verify Definition~2.2 of \cite{EN08} for $A\rtimes_{\alpha} G$. So let
			we are given a finite set $F \subseteq A\rtimes_{\alpha} G$, $\varepsilon > 0$, and
			a nonzero positive element $x \in A\rtimes_{\alpha} G$.
			We have to find a unital C*-subalgebra  $D$ of
			$A\rtimes_{\alpha} G$ with unit $e$ such that:
			
			\begin{enumerate}
				\setcounter{enumi}{\value{TmpEnumi}}
				\item\label{cor_ot_it4}
				$\|ea -ae\|<\varepsilon$ for all $a\in F$,
				
				\item\label{cor_ot_it5}
				$eFe \subseteq_{\varepsilon} D$, and
				
				\item\label{cor_ot_it6}
				$1-e \precsim x$ in $A\rtimes_{\alpha} G$.
			\end{enumerate}
			We may assume that $\varepsilon<1$ and $\|a\|\leq 1$
			for all $a\in F$. Applying Proposition~\ref{prop_arc}
			with $\varepsilon/2$ in place of $\varepsilon$, with
			1 in place of $z$ and $y$, and with $x$ and $F$ as
			given, we obtain a nonzero projection $f\in A$,
			a nonzero projection $e$ in $A\rtimes_{\alpha} G$,
			a unital C*-subalgebra $D$ of $A\rtimes_{\alpha} G$
			with unit $e$, and a $*$-isomorphism
			$\phi: M_{n} \otimes fAf \rightarrow  D$
			satisfying
			Conditions~\eqref{prop_arc_it1}--\eqref{prop_arc_it6}
			in that proposition. Note that, since $A$ is unital,
			Proposition~\ref{prop_arc} holds also without the assumption
			of Property~(SP),
			by Remark~\ref{prop_arc}.
			
			We claim that the C*-subalgebra
			$D$ and its  unit $e$ satisfy Conditions~\eqref{cor_ot_it4}--\eqref{cor_ot_it6} above.
			Note that, by the properties of the class
			$\mathcal{C}$, we have $D\in \mathcal{C}$.
			Condition~\eqref{cor_ot_it4} follows from
			the first part of Condition~\eqref{prop_arc_it3} of
			Proposition~\ref{prop_arc}.
			To see \eqref{cor_ot_it5}, let $a\in F$.
			Then by Condition~\eqref{prop_arc_it3} of
			Proposition~\ref{prop_arc} there is $b_{1}\in D$
			such that $\|ea-b_{1}\| < \varepsilon/2$,
			and that $\| ea-ae\|<\varepsilon/2$. Hence,
			\[
			\| eae- b_{1} \| \leq \| eae-ea\|+\| ea-b_{1}\|
			< \tfrac{\varepsilon}{2}+\tfrac{\varepsilon}{2}
			=\varepsilon.
			\]
			Therefore, $eFe \subseteq_{\varepsilon} D$, which is
			\eqref{cor_ot_it5}.
			Finally, Condition~\eqref{cor_ot_it6} follows
			from Condition~\eqref{prop_arc_it5} of
			Proposition~\ref{prop_arc} (see Remark~\ref{prop_arc}).
		\end{proof}
		
		The following proposition is useful to generalize results about the structure of crossed products
		of unital C*-algebras
		to the nonunital case.
		
		\begin{prop}\label{propcross}
			Let $G$ be a finite group and let $\mathcal{C}$ be a class of simple (separable) C*-algebras with
			the following properties:
			\begin{enumerate}
				\item\label{propcross_1}
				if $A$ is a simple (separable) C*-algebra and $p\in A$ is a nonzero projection, then
				$A\in\mathcal{C}$ if and only if $pAp\in\mathcal{C}$ (in particular, this is the case
				if $\mathcal{C}$ is closed under Morita equivalence);
				\item\label{propcross_2}
				if $A\in\mathcal{C}$ is unital and $\alpha$ is an action of $G$ on $A$ with the tracial
				Rokhlin property then $A\rtimes_{\alpha} G\in \mathcal{C}$;
				\item\label{propcross_3}
				if $A\in \mathcal{C}$ and $B$ is a C*-algebra with $A \cong B$, then
				$B\in \mathcal{C}$.
			\end{enumerate}
			Then $\mathcal{C}$ is closed under  taking crossed products by actions of $G$ with the
			tracial Rokhlin property (i.e., \eqref{propcross_2} above holds
			without the unital assumption).
		\end{prop}
		
		\begin{proof}
			Let  $\mathcal{C}$ be a class of simple C*-algebras as in the statement.
			Let $A\in \mathcal{C}$ and  let $\alpha \colon G \to \mathrm{Aut}(A)$
			be an action with the tracial Rokhlin property.
			We show that $A\rtimes_{\alpha} G\in \mathcal{C}$.
			We may assume that  $A$ is nonzero.
			By Lemma~\ref{lemtrp1},   there exists a nonzero projection $p$ in $A^{\alpha}$.
			Let $\beta:G\to \mathrm{Aut}(pAp)$ be the restriction of $\alpha$ to $pAp$.
			Then $\beta$ has the tracial Rokhlin property, by \cite[Proposition~4.2]{forgol}.
			By Condition~\eqref{propcross_1}, $pAp \in \mathcal{C}$. Thus,
			by Condition~\eqref{propcross_2},
			$pAp \rtimes _{\beta} G \in \mathcal{C}$. Observe that
			$pAp \rtimes _{\beta} G\cong p(A \rtimes _{\alpha} G)p$.
			In fact, the map $\varphi: pAp \rtimes _{\beta} G\to p(A \rtimes _{\alpha} G)p$
			defined by $\varphi(\sum_{g\in G}b_{g}u_{g})=\sum_{g\in G}b_{g}\delta_{g}$,
			where $b_{g}\in pAp$ for all $g\in G$, is easily seen to be a surjective $*$-isomorphism.
			Thus, by Condition~\eqref{propcross_3}, $p(A \rtimes _{\alpha} G)p \in \mathcal{C}$.
			Now Condition~\eqref{propcross_2} implies that
			$ A \rtimes _{\alpha} G \in \mathcal{C}$.
			
			If we assume that $\mathcal{C}$ is a class of \emph{separable} C*-algebras
			satisfying Conditions~\eqref{propcross_1} and \eqref{propcross_3},
			and satisfying the separable version of Condition~\eqref{propcross_2},
			then the  argument above again works, since
			$ A \rtimes _{\alpha} G$ is separable if  $A$ is separable and $G$ is finite.
		\end{proof}
		
		This following lemma is a variant of \cite[Lemma~3.6.9]{Lin.book}.
		
		\begin{lma}\label{zerodivisor}
			Let $A$ be a unital C*-algebra and suppose that $a\in A_{sa}$ is not invertible.
			Then, for any $\varepsilon>0$, there exists a nonzero  element $a' \in A_{sa}$
			such that $\|a-a'\|<\varepsilon$ and $a'$ is a zero divisor.
		\end{lma}
		
		\begin{proof}
			We assume that $0<\varepsilon<\|a\|$. Consider the following function $f:\R\rightarrow\R$:
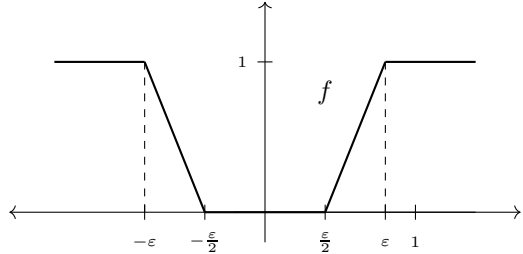
			\begin{center}
				\begin{tikzpicture}[scale=2]
					\draw[->](-0.2,0)--(1.7,0);
					\draw[->](0,-0.2)--(0,1.4);
					\draw[<-](-1.7,0)--(1.4,0);
					
					\draw(1,-0.05)--(1,0.05);
					\draw(0.8,-0.05)--(0.8,0.05);
					\draw(0.4,-0.05)--(0.4,0.05);
					\draw(-0.4,0.05)--(-0.4,-0.05);
					\draw(-0.8,-0.05)--(-0.8,0.05);
					\node at (1,-0.2){\tiny{$1$}};
					\node at (0.8,-0.2){\tiny{$\varepsilon$}};
					\node at (0.4,-0.2){\tiny{$\frac{\varepsilon}{2}$}};
					\node at (-0.4,-0.2){\tiny{$-\frac{\varepsilon}{2}$}};
					\node at (-0.8,-0.2){\tiny{$-\varepsilon$}};
					\draw(-0.05,1)--(0.05,1);
					\node at (-0.15,1){\tiny{$1$}};
					\draw[thick](0.4,0)--(0.8,1);
					\draw[thick](0.8,1)--(1.4,1);
					\draw[thick](-0.4,0)--(0.4,0);
					\draw[dashed](0.8,0)--(0.8,1);
					\node at (0.4,0.8){{$f$}};
					\draw[thick](-0.4,0)--(-0.8,1);
					\draw[thick](-0.8,1)--(-1.4,1);
					\draw[dashed](-0.8,0)--(-0.8,1);
				\end{tikzpicture}
			\end{center}

			Using functional calculus, we take $a'=f(a)a$ which is nonzero and
			satisfies $\|a-a'\|<\varepsilon$. It remains to show that $a'$ is a zero divisor.
			Let $g:\R\rightarrow\R$ be the following continuous function:
			\[
			g(t)=
			\begin{cases}
				0  & |t|>\frac{\varepsilon}{2},\\
				\frac{-2}{\varepsilon}t+1 &  0\leq t \leq \frac{\varepsilon}{2},\\
				\frac{2}{\varepsilon}t+1&  -\frac{\varepsilon}{2}\leq t <0.
			\end{cases}
			\]
			The diagram of $g$ is the following:
			
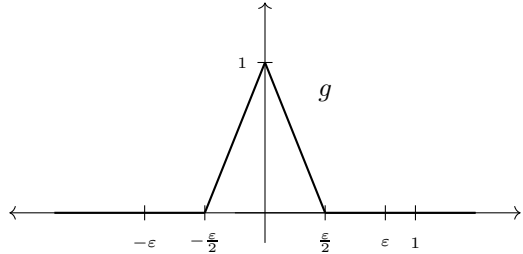
\begin{center}
				\begin{tikzpicture}[scale=2]
					\draw[->](-0.2,0)--(1.7,0);
					\draw[->](0,-0.2)--(0,1.4);
					\draw[<-](-1.7,0)--(1.4,0);
					
					\draw(1,-0.05)--(1,0.05);
					\draw(0.8,-0.05)--(0.8,0.05);
					\draw(0.4,-0.05)--(0.4,0.05);
					\draw(-0.4,0.05)--(-0.4,-0.05);
					\draw(-0.8,-0.05)--(-0.8,0.05);
					\node at (1,-0.2){\tiny{$1$}};
					\node at (0.8,-0.2){\tiny{$\varepsilon$}};
					\node at (0.4,-0.2){\tiny{$\frac{\varepsilon}{2}$}};
					\node at (-0.4,-0.2){\tiny{$-\frac{\varepsilon}{2}$}};
					\node at (-0.8,-0.2){\tiny{$-\varepsilon$}};
					\draw(-0.05,1)--(0.05,1);
					\node at (-0.15,1){\tiny{$1$}};
					\draw[thick](0,1)--(0.4,0);
					\draw[thick](0,1)--(-0.4,0);
					\draw[thick](-0.4,0)--(-1.4,0);
					\node at (0.4,0.8){{$g$}};
					\draw[thick](0.4,0)--(1.4,0);
				\end{tikzpicture}
			\end{center}
			
			We take $b=g(a)$ which is nonzero and satisfies
			$a'b=ba'=0$.
		\end{proof}
		
		In \cite{EN08}, Elliott and Niu showed that if $A\in\text{TA}\mathcal{C}$ where $\mathcal{C}$
		is a class of unital C*-algebras with stable rank one, then $A\in\mathcal{C}$.
		Also, \cite[Theorem~3.5]{Osaka} states (without proof) that this statement is valid when
		$\mathcal{C}$ is a class of
		unital C*-algebras with real rank zero but with an extra assumption that
		$A$ is stably finite. Here we remove this assumption
		and provide a complete proof using Lemma~\ref{zerodivisor}.
		
		\begin{thm}\label{Ellirr0zerodivisor}
			Let $\mathcal{C}$ be a class of unital C*-algebras with real rank zero.
			Then any simple unital C*-algebra in the class TA$\mathcal{C}$ has real rank zero.
		\end{thm}
		
		\begin{proof}
			Let $A$ be a simple unital infinite dimensional C*-algebra in TA$\mathcal{C}$.
			We show that $A$ has real rank zero.
			If $A$ does not have Property~(SP), then there is a nonzero positive element
			$a\in A$ such that the only projection of $\overline{aAa}$ is the zero projection.
			Let $x\in A_{sa}$ and $\varepsilon >0$.
			Since $A\in\text{TA}\mathcal{C}$, applying  Definition~2.2 of \cite{EN08} with $F= \{x\}$,
			and with $\varepsilon>0$ and $a$ as given, we get a nonzero projection $p\in A$ and a C*-subalgebra $D\subseteq A$
			with unit $p$ such that $D \in \mathcal{C}$, $pFp \subseteq_{\varepsilon} D$, and
			$1_{A}-p\precsim_{A} a$.
			As $\overline{aAa}$ is projectionless, the latter implies that $p=1_{A}$.
			Hence, there exists  $y\in D$  such that $\|x-y\|<\varepsilon$.
			We may assume that $y$ is self-adjoint since we can consider $\frac{y+y^*}{2}$ instead of $y$.
			Since $D$ has real rank zero, there is an invertible element $z\in D_{sa}$ such that $\|z - y\| < \varepsilon$.
			Then $\| x - z \| < 2\varepsilon$. Therefore, $A$ has real rank zero.
			
			Now we deal with the case that $A$ has  Property~(SP).
			Let $a$ be an element of $A_{sa}$ which is not invertible and let $\varepsilon > 0$.
			By Lemma~\ref{zerodivisor}, there exists a self-adjoint zero-divisor $a'$ such that
			$\|a-a'\|<\frac{\varepsilon}{2}$.
			Thus, it is enough to prove that $a'$ can be approximated by invertible self-adjoint elements of $A$.
			
			Since $A$ has Property~(SP), there is a nonzero projection $e$ which is orthogonal to $a'$.
			Indeed, let $b\in A_{sa}\setminus \{0\}$ satisfy $ba'=0$. Then  $b^2 a'=0$.
			Let $e\in\overline{b^2Ab^2}$ be a nonzero projection. Thus $ea'=0=a'e$.
			
			We may assume that both  classes $\mathcal{C}$ and TA$\mathcal{C}$ are closed under
			taking unital hereditary C*-subalgebras.
			In fact, let $\mathcal{C'}=\{pBp: B\in\mathcal{C}\ \text{and}\ p\in B \ \text{is a projection}\}$.
			Then $\mathcal{C}\subseteq\mathcal{C'}$ and $\mathcal{C'}$
			is a class of unital C*-algebras with real rank zero which is closed under taking unital hereditary C*-subalgebras,
			and hence so is  TA$\mathcal{C'}$, by   \cite[Lemma~2.3]{EN08}.
			Since $A$ is in TA$\mathcal{C}\subseteq$TA$\mathcal{C'}$,
			we may consider $\mathcal{C'}$ instead of $\mathcal{C}$.
			
			We claim that there are nonzero orthogonal projections
			$e_1 , e_2 \in eAe$ such that $e=e_1+e_2$ and $e_2\precsim_{A} e_1$. To see this,
			first    \cite[Lemma~1.10]{Ph11} implies that there exist nonzero Murray-von~Neumann equivalent  orthogonal projections
			$p_1 , p_2 \in eAe$.
			Put $e_1=e-p_2$ and $e_2=p_2$. Then $e_2 \sim p_1 \leq e-p_2=e_1$, and so  $e_2\precsim_{A} e_1$.
			Also, $e_1$ and $e_2$ are nonzero orthogonal projections with $e=e_1+e_2$.
			
			Note that $a'\in(1-e_1)A(1-e_1)$ as $a'e_1=a'ee_1=0$. Since the class TA$\mathcal{C}$ is closed under
			taking unital hereditary C*-subalgebras, $(1-e_1)A(1-e_1)$ is   in TA$\mathcal{C}$.
			So, there is a nonzero projection $p\in (1-e_1)A(1-e_1)$ and a C*-subalgebra
			${C\subseteq (1-e_1)A(1-e_1)}$ with unit $p$ such that $C\in\mathcal{C}$
			and, with $y= (1-e_1-p)a'(1-e_1-p)$, the following hold:
			\begin{enumerate}
				
				\item\label{tacequi1}
				$\|a'-(pa'p+y)\|= 2\|a'p - p a'\|<\frac{\varepsilon}{4}$;
				
				\item\label{tacequi2} $pa'p\in_{\frac{\varepsilon}{4}}C$;
				
				\item\label{tacequi3} $1-e_1-p\precsim e_2$ in $(1-e_1)A(1-e_1)$.
			\end{enumerate}
			
			Using \eqref{tacequi2} and that $C$ has real rank zero, we get an invertible self-adjoint element $b\in C = pAp$
			such that $\|b-pa'p\|<\frac{\varepsilon}{2}$.
			
			By \eqref{tacequi3} and the properties of $e_1$ and $e_2$, we have
			$1-e_1-p\precsim_{A} e_2\precsim_{A} e_1$. Hence
			there is a  partial isometry $\mu \in A$ such that $\mu\mu^*=1-e_1-p$ and $\mu^*\mu\leq e_1$.
			In particular, $\mu\in (1-e_1 - p)A e_1$. We set
			$z=y+(\frac{\varepsilon}{2})\mu+(\frac{\varepsilon}{2})\mu^*+(\frac{\varepsilon}{2})(e_1-\mu^*\mu)$, with the matrix form
			\[
			z=
			\begin{pmatrix}
				(\varepsilon/2)(e_1-\mu^*\mu) & 0 & 0\\
				0 & 0 & (\varepsilon/2)\mu^*\\
				0 & (\varepsilon/2)\mu & y
			\end{pmatrix}.
			\]
			\medskip
			
			Then $z$ is an invertible self-adjoint element of the C*-algebra $(1-p)A(1-p)$. In fact,
			the inverse of $z$ is
			$z^{-1}=(\frac{2}{\varepsilon})(e_1-\mu^*\mu) +\frac{-4}{\varepsilon^2}\mu^* y\mu +(\frac{2}{\varepsilon})\mu^*+(\frac{2}{\varepsilon})\mu$
			with the matrix form
			\[
			z^{-1}=
			\begin{pmatrix}
				(2/\varepsilon)(e_1-\mu^*\mu) & 0 & 0\\
				0 & (-4/\varepsilon^2) \mu^* y \mu & (2/\varepsilon)\mu^*\\
				0 & (2/\varepsilon)\mu & 0
			\end{pmatrix}.
			\]
			
			Note that $z+b$ is an invertible self-adjoint element of $A$, and by \eqref{tacequi1} and ${\|b-pa'p\|<\frac{\varepsilon}{2}}$
			we get
			$$	\|a'-(z+b)\| \leq\|a'-(pa'p+y)\|+\|pa'p-b\|+\|y-z\|<\frac{\varepsilon}{4}+\frac{\varepsilon}{2}+\frac{3\varepsilon}{2}.$$
			Therefore, $\|a-(z+b)\| < 3\varepsilon$.
			This shows that $A$ has real rank zero.
		\end{proof}
		
		
		\begin{thm}\label{thm_rr0}
			Let $A$ be a simple $\sigma$-unital C*-algebra with real rank zero and let $\alpha$ be an action of a finite group
			$G$ on $A$ with the tracial Rokhlin property. Then $A\rtimes_{\alpha} G$ and the fixed point algebra $A^{\alpha}$ have real rank zero.
		\end{thm}
		
		\begin{proof}
			First, we deal with the unital case.
			Suppose that $\mathcal{C}$ is the class of all simple unital C*-algebras with real rank zero.
			As the real rank zero passes to matrix algebras and hereditary C*-subalgebras \cite{BP91},
			if $A$ is a simple unital C*-algebra with real rank zero and  $\alpha$ is as in the statement,
			then by  Corollary~\ref{cor_ot}, $A\rtimes_{\alpha} G \in \mathrm{TA}\mathcal{C}$.
			Applying Theorem~\ref{Ellirr0zerodivisor}, we see that $A\rtimes_{\alpha} G \in \mathcal{C}$,
			that is, $A\rtimes_{\alpha} G$ has real rank zero. Since $A^{\alpha}$ is isomorphic to
			a corner of $A\rtimes_{\alpha} G$, it has also real rank zero.
			
			\medskip
			
			Now suppose that $\mathcal{C}$ is the class of all simple $\sigma$-unital C*-algebras with real rank zero.
			By the preceding paragraph, $\mathcal{C}$ satisfies Condition~\eqref{propcross_2} of Proposition~\ref{propcross}.
			To see Condition~\eqref{propcross_1}, let $A \in \mathcal{C}$ and
			let $p$ be a nonzero projection in $A$. As $pAp$ and $A$ are $\sigma$-unital and Morita equivalent,
			\cite[Theorem~3.8]{BP91} implies that $A$ has real rank zero if and only if $pAp$ has real rank zero.
			Therefore,  if $A$ is a simple  $\sigma$-unital C*-algebra with real rank zero and  $\alpha$ is as in the statement,
			then by  Proposition~\ref{propcross}, $A\rtimes_{\alpha} G \in \mathcal{C}$.
			Moreover, since $A\rtimes_{\alpha} G$ is simple, the fixed point algebra $A^{\alpha}$ is isomorphic to a full hereditary
			C*-subalgebra of $A\rtimes_{\alpha} G$ (by  \cite[Corollary~3.3]{forgol}). Then by
			\cite[Corollary~2.8]{BP91},  $A^{\alpha}$
			has real rank zero.	
		\end{proof}
		
		Now we deal with the preservation of stable rank one.
		The analogue of Theorem~\ref{Ellirr0zerodivisor} for stable rank one
		which is needed in the proof of the following theorem is
		\cite[Theorem~4.3]{EN08}.

		\begin{thm}\label{sr1}
			Let $A$ be a simple $\sigma$-unital C*-algebra with stable rank one and let $\alpha$ be an action of a finite group $G$ on $A$ with the tracial Rokhlin property. Then $A\rtimes_{\alpha} G$ and the fixed point algebra $A^{\alpha}$ have stable rank one.
		\end{thm}
		
		\begin{proof}
			First we consider the unital case.
			Suppose that $\mathcal{C}$ is the class of all simple unital
			C*-algebras with stable rank one. Since $\mathcal{C}$ is closed under the three conditions of Corollary~\ref{cor_ot},
			it follows that $A\rtimes_{\alpha} G\in \mathrm{TA}\mathcal{C}$
			whenever $A\in\mathcal{C}$ and
			$\alpha$ is as in the statement. By \cite[Theorem~4.3]{EN08}, $A\rtimes_{\alpha} G \in \mathcal{C}$, that is,
			it has stable rank one.
			
			Now let   $\mathcal{C}$ be the class of all simple
			$\sigma$-unital C*-algebras with stable rank one.
			We show that $\mathcal{C}$ is closed under
			Condition~\eqref{propcross_1} of Proposition~\ref{propcross}.
			Let $A \in \mathcal{C}$ and take a nonzero projection $p\in A$.
			Put $B = pAp$.
			Since $A$ and $B$ are simple and $\sigma$-unital,
			it follows that
			$A\otimes\mathcal{K}\cong B\otimes\mathcal{K}$
			\cite[Theorem~2.8]{Br77}.
			Then by \cite[Theorem~3.6]{rieff83},
			$A$ has stable rank one if and only if so does $B$.
			Condition~\eqref{propcross_2} of Proposition~\ref{propcross}
			follows from the preceding paragraph.
			Therefore, if $A$ and $\alpha$ are as in the statement
			then
			$A\rtimes_{\alpha} G$ has stable rank one.
			
			To show that $A^{\alpha}$ has stable rank one, note that
			$A\rtimes_{\alpha} G$ is simple and the fixed point algebra
			$A^{\alpha}$ is isomorphic to a full hereditary C*-subalgebra
			of $A\rtimes_{\alpha} G$ (by  \cite[Corollary~3.3]{forgol}).
			On the other hand,
			both $A\rtimes_{\alpha} G$  and $A^{\alpha}$
			are $\sigma$-unital since so is $A$. Hence
			by \cite[Theorem~2.8]{Br77},
			$(A\rtimes_{\alpha} G)\otimes\mathcal{K} \cong
			A^{\alpha}  \otimes\mathcal{K}$.
			By  \cite[Theorem~3.6]{rieff83}, $A^{\alpha}$   has stable rank one if and only if so does
			$A^{\alpha}\otimes\mathcal{K}$.
			As $A\rtimes_{\alpha} G$ has stable rank one,
			we see that  $A^{\alpha}$ has stable rank one.
		\end{proof}
		
		\section{Equivalence of tracial and weak tracial Rokhlin properties}
		
		In this section, we show that under some assumption on $A'\cap A_{\omega}$, the weak tracial Rokhlin property and the tracial Rokhlin property are equivalent when $\alpha$ is an action of a finite group $G$ on $A$.
		
		We recall some definitions concerning the \textit{Kirchberg's central sequence algebra}  $\mathrm{F}(A)$. Let $\omega$ be a fixed
		free ultrafilter on $\N$. We refer the reader to   \cite[Appendix~B]{capraro15} for details and notation about ultrafilters and ultralimits.
		For a C*-algebra $A$,  let
		\begin{align*}
			c_{\omega}(A)&:=\big\{(a_1,a_2,\dots)\in\ell^{\infty}(A) \colon\,\lim_{n\in\omega}\|a_n\|=0\big\},\\
			A_{\omega}&:=\ell^{\infty}(A)/c_{\omega}(A).
		\end{align*}
		We define the two-sided annihilator of $A$ in $A_{\omega}$ as   $$\text{Ann}(A)=\text{Ann}(A,A_{\omega}):=\{b\in A_{\omega}
		\colon bA=\{0\}=Ab\},$$ which is easy to see that is an ideal of $A'\cap A_{\omega}$. The central sequence algebra then in defined as $\text{F}(A):=(A'\cap A_{\omega})/\text{Ann}(A)$ \cite{Kirch}.
		
		If $\alpha \colon G \rightarrow\mathrm{Aut} (A)$  is an action
		of a finite group $G$ on $A$, then we denote by $\alpha _{\omega}$ the induced action of $G$ on $A_{\omega}$
		defined by $(\alpha_{\omega})_g:\ell^{\infty}(A)/c_{\omega}(A)\longrightarrow\ell^{\infty}(A)/c_{\omega}(A)$, $(a_n)_{n\in\N}+c_{\omega}(A)\mapsto(\alpha_{g}(a_n))_{n\in\N}+c_{\omega}(A)$, for all  $g\in G$.
		
		The following lemma should be known. We give proof for the convenience of the reader.
		
		\begin{lma}\label{normultra}
			Let $A$ be a  C*-algebra and let $(a_n)_{n\in\N}$ be a sequence in $\ell^{\infty}(A)$. If $\omega$ is a fixed free ultrafilter on $\N$, then  $\|(a_n)_{n\in\N}+c_{\omega}(A)\|=\lim_{n\in\omega}\|a_n\|$.
		\end{lma}
		\begin{proof}
			The idea of the proof is similar to that of \cite[Lemma~6.1.3]{RLL00}. We provide the details for completeness.
			Set $\alpha=\lim_{n\in\omega}\|a_n\|$ and $\beta=\|(a_n)+c_{\omega}\|$. For all
			$(b_n) \in c_{\omega}(A)$ we have
			\begin{align*}
				\|(a_n)+(b_n)\|&=\sup_{n\in\N}\|a_n+b_n\|\ge\lim_{n\in\omega}\|a_n+b_n\|\qquad \\
				&\ge\lim_{n\in\omega}\|a_n\|-\lim_{n\in\omega}\|b_n\|\qquad\qquad\qquad
			\end{align*}
			But $\lim_{n\in\omega}\|b_n\|=0$, and hence $\beta\ge\alpha$.
			
			Now, let $\varepsilon>0$ and consider $\mathcal{F}_{\varepsilon}:=\{n\in\N:|\alpha-\|a_n\||<\varepsilon\}\in\omega$. We define a sequence $(b_n)$  by
			\[
			b_n=
			\begin{cases}
				a_n  &n\notin\mathcal{F}_{\varepsilon},\\
				0    &n\in\mathcal{F}_{\varepsilon}.
			\end{cases}
			\]
			Then  $\lim_{n\in\omega}\|b_n\|=0$ and so $(b_n)\in c_{\omega}$.
			We have
			
			$$\beta=\|(a_n)+c_{\omega}\|\leq\|(a_n)-(b_n)\|\leq\alpha+\varepsilon.$$
			Letting $\varepsilon\rightarrow 0$,
			we get $\beta \leq \alpha$, finishing the proof.	
		\end{proof}
		
		The following is the main result of this section (Theorem~\ref{Ultrr0} of Section 1).
		
		\begin{thm}\label{ultr}
			Let $A$ be a simple separable C*-algebra and let $\alpha$ be an action
			of a finite group $G$ on $A$ with the weak tracial Rokhlin property. If $A'\cap A_{\omega}$ has
			real rank zero for some free ultrafilter  $\omega$ on $\N$, then $\alpha$ has the tracial Rokhlin property.
		\end{thm}
		\begin{proof}
			Take $\varepsilon>0$, a finite subset $F\subseteq A$, and $x,y\in A_+$ with $x\ne 0$. By the analogue of
			\cite[ Proposition~3.10]{forgol} for
			$A'\cap A_{\omega}$, there are orthogonal positive contractions $(f_g)_{g\in G}$ in $A' \cap A_{\omega}$ such that, with $f=\sum_{g \in G}f_g$, we have
			\begin{enumerate}
				\item\label{ult1}  $(\alpha_{\omega})_g(f_h)=f_{gh}$ for all $g,h\in G;$
				\item\label{ult2} $y^2-yfy\precsim_{A_{\omega}} x$;
				\item $\|fxf\|=1.$
				\setcounter{TmpEnumi}{\value{enumi}}
			\end{enumerate}
			Since $A' \cap A_{\omega}$ has real rank zero, there is a projection $p_e\in\overline{f_e(A' \cap A_{\omega})f_e}$ such that $\|p_ef_e-f_e\|<\delta$ and $\|p_ef_ep_e-f_e\|<\delta$ where $\delta=\frac{\varepsilon}{2|G|(\| y \| + 1)}$. Put $p_g=(\alpha_{\omega})_g(p_e)$ for all $g\in G$. Then $p_g\in\overline{f_g(A' \cap A_{\omega})f_g}$ and so $p_gp_h=0$ for $g,h\in G$ with $g\ne h$. Also, $(\alpha_{\omega})_g(p_h)=p_{gh}$ for all $g,h\in G$.
			In particular, by \eqref{ult1},
			$\|p_g f_g -f_g\|<\delta$ and $\|p_g f_g p_g -f_g\|<\delta$.
			We put $p=\sum_{g\in G} p_{g}$.
			We show that
			\begin{enumerate}
				\setcounter{enumi}{\value{TmpEnumi}}
				\item\label{ult4}
				$(y^2-ypy-\frac{\varepsilon}{2})_+\precsim_{A_{\omega}} x$;
				\item\label{ult5}
				$\|pxp\|>1-\varepsilon.$
				\setcounter{TmpEnumi}{\value{enumi}}
			\end{enumerate}
			For (\ref{ult4}), first we have
			\begin{align*}
				\|y^2-ypfpy-(y^2-yfy)\|
				&\leq \|y\|^{2} \cdot \|pfp-f\|\\
				&=  \|y\|^{2} \sum_{g\in G}\|p_g f_g p_g - p_g\|\\
				&=  \|y\|^{2} \sum_{g\in G}\|p_e f_e p_e - p_e\|\\
				& < \|y\|^{2}  |G|\delta =\frac{\varepsilon}{2}.
			\end{align*}
			Also, $y^2-ypy\leq y^2-ypfpy$, since $f\leq1$ implies that
			$ypfpy\leq ypy$. Hence
			$(y^2-ypy-\frac{\varepsilon}{2})_+\precsim_{A_{\omega}}
			(y^2-ypfpy-\frac{\varepsilon}{2})_+$ by \cite[Lemma~1.7]{Ph14}.
			This together with the inequality above and (\ref{ult2}) give
			$$\left(y^2-ypy-\frac{\varepsilon}{2}\right)_+\precsim_{A_{\omega}}
			\left(y^2-ypfpy-\frac{\varepsilon}{2}\right)_+\precsim_{A_{\omega}}  y^2-yfy\precsim_{A_{\omega}}  x,$$
			which is \eqref{ult4}.
			To see \eqref{ult5}, first using $\|p_gf_g-f_g\|<\delta$ we have
			$$\|f_gxf_g\|\sim_{2\delta}\|f_gp_gxp_gf_g\|\leq\|p_gxp_g\|.$$
			Thus $\|p_gxp_g\|>\|f_gxf_g\|-2\delta$. Hence
			\begin{align*}
				\|pxp\|=\max_{g\in G}\|p_gxp_g\|&>\max_{g\in G}\|f_gxf_g\|-2\delta\\
				&=\|fxf\|-2\delta=1-2\delta\\
				&\ge 1-\varepsilon.
			\end{align*}
			This is \eqref{ult5}.
			
			Now let $(p_{g,n})_{n\in\N}\in\ell^{\infty}(A)$ be a representing sequence for $p_g,\, g\in G$.
			We may assume that each $p_{g,n}$ is a self-adjoint element.
			By (\ref{ult4}), there is $v\in A_{\omega}$
			such that $\|(y^2-ypy-\frac{\varepsilon}{2})_+-vxv^*\|<\frac{\varepsilon}{2}$.
			Let $(v_n)_{n\in\N}\in\ell^{\infty}(A)$ be a representing sequence for $v$. Since $\|y^2-ypy-vxv^*\|<\varepsilon$,
			if we denote $q_n=\sum_{g \in G}p_{g,n}$, $g\in G$,
			then using Lemma~\ref{normultra}
			we get
			\begin{enumerate}
				\setcounter{enumi}{\value{TmpEnumi}}
				\item \label{aomeg6} $\lim_{n\in\omega}\|y^2-yq_ny-v_nxv_n^*\|<\varepsilon$.
				\setcounter{TmpEnumi}{\value{enumi}}
			\end{enumerate}
			Moreover, by Lemma~\ref{normultra},
			using $p_g \in A' \cap A_{\omega}$, $(\alpha_{\omega})_g (p_h)= p_{gh}$,
			\eqref{ult5}, $p_g p_h =0$ for $g\neq h$, and
			$p_{g}^2 = p_g$, respectively, we see that

			\begin{enumerate}

				\setcounter{enumi}{\value{TmpEnumi}}
				\item $\lim_{n\in\omega}\|p_{g,n}a-ap_{g,n}\|=0$, for all $g\in G$ and all $a\in A$;
				\item $\lim_{n\in\omega}\|\alpha_{g}(p_{h,n})-p_{gh,n}\|=0$ for all $g,h\in G$ with $g\ne h$;
				\item $\lim_{n\in\omega}\|q_nxq_n\|>1-\varepsilon$;
				\item $\lim_{n\in\omega}\|p_{g,n}p_{h,n}\|=0$ for all $g,h\in G$ with $g\ne h$;
				\item \label{aomeg11} $\lim_{n\in\omega}\|p_{g,n}^2-p_{g,n}\|=0$, for all $g\in G$.
				\setcounter{TmpEnumi}{\value{enumi}}
			\end{enumerate}
			Now we put
			\begin{align*}
				A_1&=\big\{n\in\N: \|y^2-yq_ny-v_nxv_n^*\|<\varepsilon\big\},\\
				A_2&=\big\{n\in\N:\|p_{g,n}a-ap_{g,n}\|<\varepsilon \ \text{for all $g\in G$ and all $a\in F$}\big\},\\
				A_3&=\big\{n\in\N: \|\alpha_{g}(p_{h,n})-p_{gh,n}\|<\varepsilon \ \text{for all $g,h\in G$} \big\},\\
				A_4&=\big\{n\in\N: \|q_nxq_n\|>1-\varepsilon\big\},\\
				A_5&=\big\{n\in\N: \|p_{g,n}p_{h,n}\|<\varepsilon\ \text{for all $g,h\in G$ with $g\neq h$}\big\},\\
				A_6&=\big\{n\in\N: \|p_{g,n}^2-p_{g,n}\|<\varepsilon\ \text{for all $g\in G$}\big\}.
			\end{align*}
			By (\ref{aomeg6})-(\ref{aomeg11}), each $A_j$ is in $\omega$ for $1\leq j\leq 6$. 
			Take $n_0\in\bigcap_{j=1}^6 A_j$, and put ${r_g=p_{g,n_0}}$. Hence $(r_g)_{g\in G}$ is a family of self-adjoint
			elements of $A$ satisfying   Remark~\ref{alttrp}. Therefore, $\alpha$ has the tracial Rokhlin property.
		\end{proof}
		
		We note that in the preceding theorem if instead of
		$A'\cap A_{\omega}$ we assume that
		$A'\cap A_{\infty}=A' \cap \left(\ell^{\infty}(A) /c_0(A)\right)$
		has real rank zero, then the same conclusion again holds.
		The proof is essentially the same.
		
		\medskip
		
		As an application of the preceding result, using
		\cite[Theorem~2.12]{Kirch}, we obtain the following.
		Note that if $A$ is unital then 
		$A'\cap A_{\omega} = \mathrm{F}(A)$ \cite[Corollary~1.10(2)]{Kirch}.
		\begin{cor}\label{cor_Kirch}
			Let $A$ be a  simple unital separable purely infinite nuclear C*-algebra, and let $\alpha$ be an action of a finite group $G$ on $A$. Then $\alpha$ has the weak tracial Rokhlin property if and only if $\alpha$ has the tracial Rokhlin property.
		\end{cor}

	\end{document}